\documentclass{amsart}
\usepackage{amssymb}

\usepackage{graphicx}
\usepackage{amscd}

\newtheorem{theorem}{Theorem}[section]
\theoremstyle{plain}
\newtheorem{lemma}[theorem]{Lemma}
\newtheorem{corollary}[theorem]{Corollary}
\newtheorem{proposition}[theorem]{Proposition}

\theoremstyle{definition}
\newtheorem{definition}[theorem]{Definition}

\theoremstyle{remark}

\numberwithin{equation}{section}

\input{tcilatex}

\begin{document}
\title[Global Left Loop Structures on Spheres]{Global Left Loop Structures on Spheres}
\author{Michael K. Kinyon}
\address{Department of Mathematics \& Computer Science\\
Indiana University South Bend\\
South Bend, IN 46634}
\email{mkinyon@iusb.edu}
\urladdr{http://www.iusb.edu/\symbol{126}mkinyon}
\date{\today}
\thanks{This paper is in final form and no version of it will be submitted for
publication elsewhere.}

\begin{abstract}
On the unit sphere $\mathbb{S}$ in a real Hilbert space $\mathbf{H}$, we
derive a binary operation $\odot $ such that $(\mathbb{S},\odot )$ is a
power-associative Kikkawa left loop with two-sided identity $\mathbf{e}_{0}$%
, i.e., it has the left inverse, automorphic inverse, and $A_{l}$
properties. The operation $\odot $ is compatible with the symmetric space
structure of $\mathbb{S}$. $(\mathbb{S},\odot )$ is not a loop, and the
right translations which fail to be injective are easily characterized. $(%
\mathbb{S},\odot )$ satisfies the left power alternative and left Bol
identities ``almost everywhere'' but not everywhere. Left translations are
everywhere analytic; right translations are analytic except at $-\mathbf{e}%
_{0}$ where they have a nonremovable discontinuity. The orthogonal group $O(%
\mathbf{H})$ is a semidirect product of $(\mathbb{S},\odot )$ with its
automorphism group. The left loop structure of $(\mathbb{S},\odot )$ gives
some insight into spherical geometry.
\end{abstract}

\maketitle

\section{Introduction}

Let $\mathbf{H}$ be a real Hilbert space with inner product $\langle \cdot
,\cdot \rangle $. Let 
\begin{equation}
\mathbb{S}=\left\{ \mathbf{x}\in \mathbf{H}:\langle \mathbf{x},\mathbf{x}%
\rangle =1\right\}  \label{sphere}
\end{equation}
be the unit sphere in $\mathbf{H}$. For $\mathbf{x},\mathbf{y}\in \mathbb{S}$%
, set 
\begin{equation}
\mathbf{x}\ast \mathbf{y}=2\langle \mathbf{x},\mathbf{y}\rangle \mathbf{x}-%
\mathbf{y}.  \label{eq:sphere_symm}
\end{equation}
(\cite{loos}, p.66). If $\mathbf{x}$ and $\mathbf{y}$ span a plane $\Pi $
through the origin, then $\mathbf{x}\ast \mathbf{y}$ is a point lying in $%
\mathbb{S}\cap \Pi $ which is obtained by reflecting $\mathbf{y}$ in $\Pi $
about the line passing through $\mathbf{x}$ and the origin. Equivalently,
thinking of $\mathbb{S}\cap \Pi $ as the great circle passing through $%
\mathbf{x}$ and $\mathbf{y}$, $\mathbf{x}\ast \mathbf{y}$ is the point on $%
\mathbb{S}\cap \Pi $ whose angle with $\mathbf{x}$ is the same as the angle
between $\mathbf{x}$ and $\mathbf{y}$, but with the opposite orientation.
The magma $(\mathbb{S},\ast )$ has the structure of a symmetric space.

\begin{definition}
\label{defn:symm-space}A \emph{symmetric space} $(M,\ast )$ is a topological
space $M$ together with a continuous binary operation $\ast :M\times
M\rightarrow M$ such that the following properties hold. For $x,y,z\in M$,

\begin{enumerate}
\item[(L1)]  (idempotence) $\ \ x\ast x=x$

\item[(L2)]  (left keyes identity) $\ \ x\ast (x\ast y)=y$

\item[(L3)]  (left distributivity) $\ \ \ x\ast (y\ast z)=(x\ast y)\ast
(x\ast z)$

\item[(L4)]  Each $a\in M$ has a neighborhood $U\subseteq M$ such that $%
a\ast x=x$ implies $x=a$ for all $x\in U$.
\end{enumerate}
\end{definition}

We make the implicit assumption that the binary operation $\ast :M\times
M\rightarrow M$ possesses as much smoothness as the space $M$ allows. Thus
if $M$ is a smooth manifold, $\ast $ is assumed to be smooth, and not just
continuous. In the smooth case, Definition \ref{defn:symm-space} is due to
Loos \cite{loos}. It is equivalent to the more standard ones (cf. \cite
{helgason}), but has the advantage of depending only on the topology. This
motivates our trivial adaptation of his definition.

The magma $(\mathbb{S},\ast )$ defined by (\ref{eq:sphere_symm}) is easily
seen to satisfy (L1), (L2) and (L3). Also, $\mathbf{x}\ast \mathbf{y}=%
\mathbf{y}$ if and only if $\mathbf{y}=\pm \mathbf{x}$, and thus (L4) holds.

Observe that the first three axioms of a symmetric space are purely
algebraic. For $x\in M$, we define a mapping $s_{x}:M\rightarrow M$, called
the \emph{symmetry about} $x$, by 
\begin{equation}
s_{x}(y)=x\ast y.  \label{symmetry}
\end{equation}
In terms of the symmetries, the keyes identity (L2) reads $s_{x}^{2}=I$,
where $I$ is the identity mapping on $M$. Therefore $(M,\ast )$ is a left
quasigroup, i.e., each symmetry $s_{x}$ is a bijection. The topological
axiom (L4) asserts that each $a\in M$ is an isolated fixed point of the
symmetry $s_{a}$.

As another example of a symmetric space, let $G$ be a topological group with
the property that $1\in G$ has a neighborhood $U\subseteq G$ such that no
element of $U$ has order 2. (For instance, if $G$ is a Lie group, then the
existence of such a neighborhood is implied by the existence of the
exponential map.) Define 
\begin{equation}
x\ast y=xy^{-1}x  \label{group-symm}
\end{equation}
for $x,y\in G$. Axioms (L1), (L2) and (L3) are easily checked. Now $1\ast
x=x^{-1}$ for all $x\in G$, and thus the assumed property of $G$ implies
that (L4) holds for $a=1$. In addition, we have $x(y\ast z)=xy\ast xz$ for $%
x,y,z\in G$ (\cite{loos}, p.65). This implies that (L4) holds for all $a\in
G $.

The case where $G$ is a (multiplicative) abelian group is of particular
interest here; in this case the symmetric space operation is simply $x\ast
y=x^{2}y^{-1}$ for $x,y\in G$. Now assume $G$ is uniquely $2$-divisible,
i.e., that the squaring map is bijective. For $x\in G$, let $x^{1/2}$ denote
the unique element of $G$ whose square is $x$. Then one can recover the
group product from the symmetric space operation by defining 
\begin{equation}
x\cdot y=x^{1/2}\ast (1\ast y)  \label{eq:abelian-op}
\end{equation}
for all $x,y\in G$. Thus if $G$ is uniquely $2$-divisible, $(G,\cdot )$ and $%
(G,\ast )$ are isotopic. It is clear from (\ref{eq:abelian-op}) that the two
essential ingredients for the isotopy are: a distinguished point $1\in G$
and a well-defined notion of ``$x^{1/2}$'' for all $x\in G$. It is useful to
note that the latter ingredient can be described entirely in terms of the
symmetric space operation. If $G$ is a uniquely $2$-divisible abelian group,
then \emph{for} $x\in G$\emph{,} \emph{there exists a unique} $x^{1/2}\in G$ 
\emph{such that} $x^{1/2}\ast 1=x$.

This idea can be generalized to a wider class of symmetric spaces.

\begin{lemma}
\label{lem:symm_quasi}Let $(M,\ast )$ be a left keyesian, left distributive,
left quasigroup, and fix $e\in M$. The following properties are equivalent.

\begin{enumerate}
\item  For each $x\in M$, there exists a unique $x^{1/2}\in M$ such that $%
x^{1/2}\ast e=x$.

\item  $(M,\ast )$ is a right quasigroup.
\end{enumerate}
\end{lemma}

\begin{proof}
(1)$\Longrightarrow $(2): Let $a,b\in G$ be given, and assume $x\in M$
satisfies $x\ast a=b$. Then using (L3) and (L2), 
\begin{equation*}
a^{1/2}\ast b=a^{1/2}\ast (x\ast a)=(a^{1/2}\ast x)\ast (a^{1/2}\ast
a)=(a^{1/2}\ast x)\ast e.
\end{equation*}
By the uniqueness assumption, $a^{1/2}\ast x=(a^{1/2}\ast b)^{1/2}$, and
thus by (L2), $x=a^{1/2}\ast (a^{1/2}\ast b)^{1/2}$. Conversely, set $%
x=a^{1/2}\ast (a^{1/2}\ast b)^{1/2}$. Then using (L3) and (L2), 
\begin{eqnarray*}
x\ast a &=&\left( a^{1/2}\ast \left( a^{1/2}\ast b\right) ^{1/2}\right) \ast
\left( a^{1/2}\ast e\right) \\
&=&a^{1/2}\ast \left( \left( a^{1/2}\ast b\right) ^{1/2}\ast e\right) \\
&=&a^{1/2}\ast \left( a^{1/2}\ast b\right) \\
&=&b.
\end{eqnarray*}
Therefore $(M,\ast )$ is a right quasigroup and hence a quasigroup.

\noindent (2)$\Longrightarrow $(1): Obvious.
\end{proof}

\begin{lemma}
\label{lem:sss_redundant}Let $(M,\ast )$ be a left distributive quasigroup.
Then $(M,\ast )$ is idempotent. Moreover, for each $a\in M$, $a\ast x=x$
implies $x=a$ for all $x\in M$.
\end{lemma}

\begin{proof}
For all $x,a\in M$, $(x\ast x)\ast a=x\ast (x\ast (x\backslash a))=x\ast a$.
Cancelling $a$, we have $x\ast x=x$, and thus $(M,\ast )$ is idempotent. If $%
a\ast x=x=x\ast x$, then clearly $x=a$.
\end{proof}

These considerations lead us to the following definition.

\begin{definition}
\label{defn:sss}A \emph{reflection quasigroup} $(M,\ast )$ is a left
keyesian, left distributive quasigroup.
\end{definition}

The term ``reflection quasigroup'' is due to Kikkawa \cite{kikkawa1}.
Reflection quasigroups are also known as ``left-sided quasigroups'',
following a convention of Robinson \cite{rob-right}. Reflection quasigroups
are also equivalent to the recently studied ``point reflection structures''
of Gabrieli and Karzel \cite{gk1} \cite{gk2} \cite{gk3}. Given a nonempty
set $P$ and a mapping $\symbol{126}:P\rightarrow \mathrm{Sym}%
(P):x\longmapsto \tilde{x}$, Gabrieli and Karzel call the pair $(P,\symbol{%
126})$ a ``point-reflection structure'' if the following hold: (i) $\forall
x,y\in P$, $(\tilde{x}\circ \tilde{x})(y)=y$; (ii) $\forall x,y\in P$, $%
\tilde{x}(y)=y$ implies $y=x$; (iii) $\forall a,b\in P$, $\exists !x\in P$
such that $\tilde{x}(a)=b$; (iv) $\forall a,b\in P$, $\exists c\in P$ such
that $\tilde{a}\circ \tilde{b}\circ \tilde{a}=\tilde{c}$. This is simply the
notion of a reflection quasigroup in different notation. Indeed, for $x,y\in
P$, set $x\ast y=\tilde{x}(y)$. Then $(P,\ast )$ is left keyesian by (i),
idempotent by (ii), and is a right quasigroup by (iii). Applying both sides
of (iv) to $a\ast b$, and using the left keyesian identity and idempotence,
we have $a\ast b=c\ast (a\ast b)$. By idempotence and the fact that $(P,\ast
)$ is a right quasigroup, $c=a\ast b$. Thus applying (iv) to $a\ast x$, and
using the left keyesian identity, we have $a\ast (b\ast x)=(a\ast b)\ast
(a\ast x)$ for all $a,b,x\in P$. This is left distributivity. Conversely, if 
$(P,\ast )$ is a reflection quasigroup, then for $x,y\in P$, set $\tilde{x}%
(y)=x\ast y$. Properties (i), (ii), and (iii) are clear. For all $a,b,x\in P$%
, we have $(\tilde{a}\circ \tilde{b}\circ \tilde{a})(x)=a\ast (b\ast (a\ast
x))=(a\ast b)\ast x$, using the left distributive and left keyesian
properties. Thus (iv) holds with $c=a\ast b$.

Lemma \ref{lem:sss_redundant} implies that a topological reflection
quasigroup is necessarily a symmetric space. However, the symmetric space $(%
\mathbb{S},\ast )$ is not a reflection quasigroup. For instance, given $%
\mathbf{a}\in \mathbb{S}$, the equation $\mathbf{x}\ast \mathbf{a}=\mathbf{a}
$ has two solutions $\mathbf{x}=\pm \mathbf{a}$. The obstructions to $(%
\mathbb{S},\ast )$ being a reflection quasigroup can also be seen in
geometric terms. If $\dim \mathbf{H}>2$, then for a fixed $\mathbf{e}\in 
\mathbb{S}$, there is a nontrivial family of great circles on $\mathbb{S}$
which connect $\mathbf{e}$ to its antipode $-\mathbf{e}$. Thus consider the
equation 
\begin{equation}
\mathbf{x}\ast \mathbf{e}=-\mathbf{e}  \label{eq:x*e=-e}
\end{equation}
in the symmetric space $(\mathbb{S},\ast )$. We have $-\mathbf{e}=2\langle 
\mathbf{e},\mathbf{x}\rangle \mathbf{x}-\mathbf{e}$, or $\langle \mathbf{e},%
\mathbf{x}\rangle \mathbf{x}=\mathbf{0}$. Hence (\ref{eq:x*e=-e}) is
equivalent to 
\begin{equation}
\langle \mathbf{e},\mathbf{x}\rangle =0.  \label{eq:<e,x>=0}
\end{equation}
Therefore any point $\mathbf{x}\in \mathbb{S}$ which is orthogonal to $%
\mathbf{e}$ will satisfy (\ref{eq:x*e=-e}). But in general, there is no
distinguished solution to (\ref{eq:x*e=-e}). (This nonuniqueness problem
holds even for the unit circle $S^{1}\subset \mathbb{R}^{2}$. However, since
there are only two solutions, a choice based on orientation leads to the
usual group structure of $S^{1}$.) Despite these concerns, reflection
quasigroups form a model for our later discussion of $(\mathbb{S},\ast )$.

The aforementioned obstructions extend to other compact symmetric spaces.
Indeed, in a smooth symmetric space, each symmetry $s_{x}$ is locally a
geodesic symmetry; i.e., for $y$ near $x$, $s_{x}(y)$ is the image of $y$
under the reflection about $x$ along a (maximal) geodesic connecting $x$ and 
$y$. For a distinguished point $e$ and any other point $x$ in a smooth
reflection quasigroup (such as a Riemannian symmetric space of noncompact
type), there is only one geodesic connecting $e$ and $x$, and thus there is
a unique point $x^{1/2}$ on that geodesic such that $x^{1/2}\ast
e=s_{x^{1/2}}(e)=x$. However, in compact symmetric spaces, the distinguished
point $e$ (or any point for that matter) has a nonempty cut locus, i.e., set
of points conjugate to $e$ \cite{klingenberg}. For each such conjugate point 
$x$, there is a nontrivial family of distinct geodesics connecting $e$ and $%
x $, and thus there is a nontrivial family of distinct points $z$ such that $%
z\ast e=x$. Since there is generally no distinguished choice for ``$x^{1/2}$%
'', compact smooth symmetric spaces are generally not reflection quasigroups.

Let $(M,\ast )$ be a reflection quasigroup, fix a distinguished element $%
e\in M$, and as before, for each $x\in M$, let $x^{1/2}=x/e$, where $/$
denotes right division in $(M,\ast )$. (Since a reflection quasigroup is
left keyesian, left division $\backslash $ agrees with the multiplication $%
\ast $.) Define 
\begin{equation}
x\cdot y=x^{1/2}\ast (e\ast y)  \label{eq:ssymm_op}
\end{equation}
for $x,y\in M$. Then $(M,\cdot )$ is a loop with identity element $e$, and
in fact, $(M,\cdot )$ is the principal $e,e$-isotope of $(M,\ast )$ \cite
{bruck} \cite{pflug}: $x\cdot y=(x/e)\ast (e\backslash y)$. Following some
definitions, we will identify the class of loops to which $(M,\cdot )$
belongs.

In a left loop, denote the left translations by $L_{x}:y\longmapsto x\cdot y$%
, and the left inner mappings by $L(x,y)=L_{x\cdot y}^{-1}L_{x}L_{y}$. A
left loop is said to have the \emph{left inverse property} if 
\begin{equation}
L_{x}^{-1}=L_{x^{-1}}  \tag{LIP}
\end{equation}
for all $x$, where $x^{-1}$ is the right (and necessarily, two-sided)
inverse of $x$. A left loop is said to have the $A_{l}$ \emph{property} if
every left inner mapping is an automorphism. A left loop with two-sided
inverses is said to have the \emph{automorphic inverse property} if 
\begin{equation}
(x\cdot y)^{-1}=x^{-1}\cdot y^{-1}  \tag{AIP}
\end{equation}
for all $x,y$. An $A_{l}$, LIP, AIP left loop is called a \emph{Kikkawa}
left loop \cite{kiechle}. A loop is called a (left) \emph{Bol loop} if it
satisfies the (left) \emph{Bol identity}: 
\begin{equation}
L_{x}L_{y}L_{x}=L_{x\cdot (y\cdot x)}  \tag{Bol}
\end{equation}
for all $x,y,z$. Bol loops necessarily satisfy LIP (see, e.g., \cite{kiechle}%
, 6.4). A Bol loop with AIP is called a \emph{Bruck loop}. Bruck loops are
necessarily Kikkawa loops$\ $(see, e.g., \cite{kiechle}, 6.6(3)). A Bruck
loop is called a \emph{B-loop} if the mapping $x\longmapsto x\cdot x$ is a
permutation.

Robinson's original definition of ``Bruck loop'' \cite{rob-diss} \cite
{rob-right} included the property that squaring is a permutation, and is
thus equivalent to what we call a B-loop. Contemporary usage of the term
``Bruck loop'' in the literature (with some exceptions) tends to be as we
have given it here, cf. \cite{kepka} \cite{kreuzer}. This usage seems to
stem from a remark of Glauberman (\cite{glauberman}, p.376).

Bruck loops are equivalent to the class of loops Ungar dubbed ``gyrogroups'' 
\cite{ungar-gyrogroup} and later ``gyrocommutative gyrogroups'' \cite
{ungar-axioms}. Bruck loops are also equivalent to ``K-loops'', which are
the additive loops of near-domains. Near-domains were introduced by Karzel 
\cite{karzel}, and the additive loop structure was later axiomatized and
named in unpublished work of Kerby and Wefelscheid (the first appearance of
the term ``K-loop'' in the literature was in a paper of Ungar \cite
{ungar-kloop}). The aforementioned equivalences have been established
independently by various authors. Kreuzer showed the equivalence of Bruck
loops with K-loops \cite{kreuzer}, and Sabinin \textit{et al} showed the
equivalence of Bruck loops with (gyrocommutative) gyrogroups \cite{sabsab}.
(The direct equivalence of gyrocommutative gyrogroups with K-loops is a
well-known folk result.)

The term ``B-loop'' was introduced by Glauberman \cite{glauberman} to
describe a finite, odd order Bol loop with the automorphic inverse property.
Having odd order implies that the squaring map is a permutation, and
Glauberman noted that for some results, it is this latter property which is
essential, not the finiteness (\cite{glauberman}, p. 374). Since the
contemporary usage of ``Bruck loop'' no longer implies that squaring is a
bijection, it is quite natural to extend Glauberman's terminology to the not
necessarily finite case. (A less acronymic alternative would be ``uniquely $%
2 $-divisible Bruck loop''.)

\begin{proposition}
\label{prop:sss=B-loop}

\begin{enumerate}
\item  Let $(M,\ast )$ be a reflection quasigroup with distinguished element 
$e\in M$, and define $x\cdot y=x^{1/2}\ast (e\ast y)$ for $x,y\in M$. Then $%
(M,\cdot )$ is a B-loop.

\item  Let $(M,\cdot )$ be a B-loop and define $x\ast y=x^{2}\cdot y^{-1}$
for $x,y\in M$. Then $(M,\ast )$ is a reflection quasigroup.
\end{enumerate}
\end{proposition}

In short, every reflection quasigroup is isotopic to a B-loop and conversely.

With different terminology and motivation than that employed here,
Proposition \ref{prop:sss=B-loop} seems to have been first obtained by D.
Robinson in his 1964 dissertation \cite{rob-diss}; the result was not
published, however, until 1979 \cite{rob-right}. Again with different jargon
and with different axioms for what we are calling a B-loop, Kikkawa
discovered the result independently in 1973 \cite{kikkawa1}. Kikkawa's
version emphasized the connection with symmetric spaces. The isotopic
relationship between reflection quasigroups and B-loops is quite natural,
and has been rediscovered other times in the literature, e.g., \cite{gk1} 
\cite{gk2} \cite{gk3}.

We have already noted that a smooth compact symmetric space $(M,\ast )$ is
not a reflection quasigroup, and thus the symmetric space operation $\ast $
is not obtained from a B-loop operation by the previously described isotopy.
However, for smooth symmetric spaces, whether compact or not, there is a
well-developed \emph{local} theory which guarantees that in a neighborhood
of any distinguished point $e$, there exists a local B-loop structure with $%
e $ as its identity element. This theory was worked out in detail primarily
by Sabinin \cite{sabinin81}; expositions can be found in \cite{mik-sab} and 
\cite{sab-book}. For a smooth reflection quasigroup, the globally-defined
B-loop operation agrees with the locally-defined operation guaranteed by the
general theory wherever the latter operation is defined.

Every Bol loop, and hence every B-loop, is left alternative, and thus the
relationship between a reflection quasigroup $(M,\ast )$ and the B-loop $%
(M,\cdot )$ can be written 
\begin{equation*}
x\ast y=x\cdot (x\cdot y^{-1})
\end{equation*}
for $x,y\in M$. In the next section, we derive a \emph{globally}-defined
binary operation $\odot :\mathbb{S}\times \mathbb{S}\rightarrow \mathbb{S}$
which is compatible with the symmetric space structure $(\mathbb{S},\ast )$
in the same sense: 
\begin{equation}
\mathbf{x}\ast \mathbf{y}=\mathbf{x}\odot (\mathbf{x}\odot \mathbf{y}^{-1})
\label{eq:compatible}
\end{equation}
for all $\mathbf{x},\mathbf{y}\in \mathbb{S}$. Our operation $\odot $ will
agree with the local B-loop operation guaranteed by general theory wherever
the latter operation is defined. However, our approach will be elementary,
in that we will not use any tools from differential geometry, nor even the
intrinsic spherical distance on $\mathbb{S}$. Rather, we will simply use (%
\ref{eq:sphere_symm}) and the Hilbert space structure.

\section{Notation}

We now introduce some notation which will be used throughout the paper. As
in \S 1, let $\mathbf{H}$ be a real Hilbert space with inner product $%
\langle \cdot ,\cdot \rangle :\mathbf{H}\times \mathbf{H}\rightarrow \mathbb{%
R}$ and corresponding norm $\left\| \mathbf{x}\right\| =\langle \mathbf{x},%
\mathbf{x}\rangle ^{1/2}$ for $\mathbf{x}\in \mathbf{H}$. For any subset $%
\mathbf{W}\subseteq \mathbf{H}$, we will denote by $\mathbf{W}^{\perp }$ the
orthogonal complement of $\mathbf{W}$ in $\mathbf{H}$. Let $\mathbf{e}_{0}$
be a distinguished unit vector in $\mathbf{H}$, and let $\mathbf{V}=\mathbf{e%
}_{0}^{\perp }$ be the orthogonal complement in $\mathbf{H}$ of $\mathbf{e}%
_{0}$. Elements of the orthogonal direct sum $\mathbf{H}=\mathbb{R}\mathbf{e}%
_{0}\oplus \mathbf{V}$ will be denoted by $\mathbf{x}=x_{0}\mathbf{e}_{0}+%
\mathbf{x}_{\perp }$ where $x_{0}=\langle \mathbf{x},\mathbf{e}_{0}\rangle $%
. In particular, for $\mathbf{V}=\mathbb{R}^{n}$, it is useful to identify $%
\mathbf{e}_{0}\in \mathbb{R}^{n+1}$ with the standard basis element $%
(1,0,\ldots ,0)^{T}$.

For any closed subspace $\mathbf{E}\subseteq \mathbf{H}$, let $\mathcal{B}(%
\mathbf{E})$ denote the set of all bounded linear operators on $\mathbf{E}$.
Given $A\in \mathcal{B}(\mathbf{E})$, the \emph{transpose} of $A$ is the
operator $A^{T}\in \mathcal{B}(\mathbf{E})$ defined by 
\begin{equation}
\langle A^{T}\mathbf{x},\mathbf{y}\rangle =\langle \mathbf{x},A\mathbf{y}%
\rangle  \label{eq:transpose}
\end{equation}
for $\mathbf{x},\mathbf{y}\in \mathbf{E}$. A bounded linear operator $A$ is 
\emph{symmetric} if $A^{T}=A$, and is \emph{orthogonal} if $A^{T}A=I$, the
identity transformation.

For a given $\mathbf{x}\in \mathbf{H}$, we define a linear functional $%
\mathbf{x}^{T}:\mathbf{H}\rightarrow \mathbb{R}$ by 
\begin{equation}
\mathbf{x}^{T}\mathbf{y}=\langle \mathbf{x},\mathbf{y}\rangle
\label{eq:functional}
\end{equation}
for $\mathbf{y}\in \mathbf{H}$. For $\mathbf{a},\mathbf{b}\in \mathbf{H}$,
we define an operator $\mathbf{ab}^{T}\in \mathcal{B}(\mathbf{H})$ by 
\begin{equation}
\mathbf{ab}^{T}\mathbf{x}=\langle \mathbf{b,x}\rangle \mathbf{a}
\label{eq:ab}
\end{equation}
for $\mathbf{x}\in \mathbf{H}$. This operator satisfies $(\mathbf{ab}%
^{T})^{T}=\mathbf{ba}^{T}$. We denote the orthogonal projection onto the
subspace $\mathbb{R}\cdot \mathbf{a}$ by 
\begin{equation}
P_{\mathbf{a}}=\frac{\mathbf{aa}^{T}}{\left\| \mathbf{a}\right\| ^{2}}
\label{eq:projection}
\end{equation}
The operator $P_{\mathbf{a}}$ is symmetric and $P_{\mathbf{a}}^{2}=P_{%
\mathbf{a}}$.

For any closed subspace $\mathbf{E}\subseteq \mathbf{H}$, let 
\begin{equation}
O(\mathbf{E})=\{A\in \mathcal{B}(\mathbf{E}):A^{T}A=I\}  \label{eq:O(E)}
\end{equation}
be the orthogonal group of $\mathbf{E}$. We will identify $O(\mathbf{V})$
with a particular subgroup of $O(\mathbf{H})$, namely 
\begin{equation}
O(\mathbf{V})\cong \left\{ A\in O(\mathbf{H}):A\mathbf{e}_{0}=\mathbf{e}%
_{0}\right\} .  \label{eq:O(V)}
\end{equation}
For any $\mathbf{a}\in \mathbf{H}$, $A\in O(\mathbf{H})$, we have 
\begin{equation}
AP_{\mathbf{a}}=P_{A\mathbf{a}}A^{T}.  \label{eq:APa=PAaA}
\end{equation}

Finally, let $J\in \mathcal{B}(\mathbf{H})$ be defined by 
\begin{equation}
J\mathbf{x}=x_{0}\mathbf{e}_{0}-\mathbf{x}_{\perp }  \label{eq:J}
\end{equation}
for $\mathbf{x}\in \mathbf{H}$. Then $J$ is symmetric and orthogonal, and
the following properties hold: 
\begin{eqnarray}
J^{2} &=&I  \label{eq:JJ=I} \\
J\mathbf{e}_{0} &=&\mathbf{e}_{0}  \label{eq:Je0=e0} \\
JV &=&V  \label{eq:JV=V} \\
I+J &=&2\mathbf{e}_{0}\mathbf{e}_{0}^{T}  \label{eq:I+J=2e0e0}
\end{eqnarray}

\section{Derivation of the Operation}

We now derive the operation $\odot $ on $\mathbb{S}$, which will turn out to
give $(\mathbb{S},\odot )$ a left loop structure. The distinguished point $%
\mathbf{e}_{0}\in \mathbb{S}$ will be the identity element. We will use (\ref
{eq:ssymm_op}) as a guide. For $\mathbf{x}\in \mathbb{S}$, we wish to define 
$\mathbf{x}^{1/2}$ to be an element of $\mathbb{S}$ satisfying $\mathbf{x}%
^{1/2}\ast \mathbf{e}_{0}=\mathbf{x}$, provided such an element exists, and
provided there is some canonical choice among such elements. For $\mathbf{u},%
\mathbf{x}\in \mathbb{S}$, we have $\mathbf{u}\ast \mathbf{e}_{0}=\mathbf{x}$
if and only if 
\begin{equation}
2\langle \mathbf{u},\mathbf{e}_{0}\rangle \mathbf{u}=\mathbf{e}_{0}+\mathbf{x%
}.  \label{eq:root_comp}
\end{equation}
If $\mathbf{x}=-\mathbf{e}_{0}$, then $\langle \mathbf{u},\mathbf{e}%
_{0}\rangle =0$, and hence $\mathbf{u}$ could be any element of $\mathbb{S}%
\cap \mathbf{V}$. Thus there is no unique choice of square root for the
antipode $-\mathbf{e}_{0}\in \mathbb{S}$. Therefore assume for now that $%
\mathbf{x}\neq -\mathbf{e}_{0}$. Taking norms of both sides of (\ref
{eq:root_comp}), we have $2\left| \mathbf{u}^{T}\mathbf{e}_{0}\right|
=\left\| \mathbf{e}_{0}+\mathbf{x}\right\| $. Therefore (\ref{eq:root_comp})
implies 
\begin{equation}
\mathbf{u}=\pm \frac{\mathbf{e}_{0}+\mathbf{x}}{\left\| \mathbf{e}_{0}+%
\mathbf{x}\right\| }.  \label{eq:u-comp}
\end{equation}
Geometrically, it is clear from (\ref{eq:u-comp}) that the preferred choice
of square root is 
\begin{equation}
\mathbf{x}^{1/2}=\frac{\mathbf{e}_{0}+\mathbf{x}}{\left\| \mathbf{e}_{0}+%
\mathbf{x}\right\| }  \label{eq:root}
\end{equation}
for $\mathbf{x}\neq -\mathbf{e}_{0}$. This can also be seen by a continuity
argument: as $\mathbf{x}\rightarrow \mathbf{e}_{0}$, we have $\mathbf{u}%
\rightarrow \pm \mathbf{e}_{0}$. Taking $\mathbf{e}_{0}^{1/2}=\mathbf{e}_{0}$
and assuming $\mathbf{x}\longmapsto \mathbf{x}^{1/2}$ to be continuous at $%
\mathbf{e}_{0}$ gives the plus sign convention in (\ref{eq:u-comp}).

With (\ref{eq:ssymm_op}) as motivation, for $\mathbf{x},\mathbf{y}\in 
\mathbb{S}$ with $\mathbf{x}\neq -\mathbf{e}_{0}$, we set 
\begin{equation}
\mathbf{x}\odot \mathbf{y}=\mathbf{x}^{1/2}\ast (\mathbf{e}_{0}\ast \mathbf{y%
}).  \label{eq:op_tmp1}
\end{equation}
Now for $\mathbf{y}\in \mathbb{S}$, we have 
\begin{equation*}
\mathbf{e}_{0}\ast \mathbf{y}=2\langle \mathbf{e}_{0},\mathbf{y}\rangle 
\mathbf{e}_{0}-\mathbf{y}=2y_{0}\mathbf{e}_{0}-\mathbf{y}=y_{0}\mathbf{e}%
_{0}-\mathbf{y}_{\perp }=J\mathbf{y}.
\end{equation*}
Therefore for $\mathbf{x},\mathbf{y}\in \mathbb{S}$ with $\mathbf{x}\neq -%
\mathbf{e}_{0}$, we have the following expression for the operation $\odot $%
: 
\begin{equation}
\mathbf{x}\odot \mathbf{y}=\left( 2P_{\mathbf{x}^{1/2}}-I\right) J\mathbf{y}.
\label{eq:op1}
\end{equation}

What remains is to find a suitable definition of $-\mathbf{e}_{0}\odot 
\mathbf{y}$. As noted, (\ref{eq:op_tmp1}) is not appropriate here, because
there is no canonical choice for ``$\left( -\mathbf{e}_{0}\right) ^{1/2}$''.
However, the correct definition of $-\mathbf{e}_{0}\odot \mathbf{y}$ will be
clear from the following remarks. Let $\Pi \subset \mathbf{H}$ be a plane
passing through $\mathbf{e}_{0}$ and the origin, and fix $\mathbf{e}_{1}\in 
\mathbb{S}\cap \Pi \cap \mathbf{V}$. For $\mathbf{x}\in \mathbb{S}\cap \Pi $%
, we have $\mathbf{x}=x_{0}\mathbf{e}_{0}+x_{1}\mathbf{e}_{1}$, where $x_{1}=%
\mathbf{x}_{\perp }^{T}\mathbf{e}_{1}$. In particular, for $\mathbf{x},%
\mathbf{y}\in \Pi $, $\mathbf{e}_{0}+\mathbf{x}=(1+x_{0})\mathbf{e}_{0}+x_{1}%
\mathbf{e}_{1}$, $\left\| \mathbf{e}_{0}+\mathbf{x}\right\|
^{2}=(1+x_{0})^{2}+x_{1}^{2}=2(1+x_{0})$, and $J\mathbf{y}=y_{0}\mathbf{e}%
_{0}-y_{1}\mathbf{e}_{1}$. Thus for $\mathbf{x},\mathbf{y}\in \mathbb{S}\cap
\Pi $, it is easy to show with a few tedious calculations using (\ref
{eq:root}) that (\ref{eq:op1}) simplifies as follows: 
\begin{equation}
\mathbf{x}\odot \mathbf{y}=(x_{0}y_{0}-x_{1}y_{1})\mathbf{e}%
_{0}+(x_{0}y_{1}+x_{1}y_{0})\mathbf{e}_{1}.  \label{eq:circle}
\end{equation}
Therefore the set $\mathbb{S}\cap \Pi $ is closed under the partial
operation $\odot $. Now the limit of the expression in (\ref{eq:circle}) as $%
\mathbf{x}\rightarrow -\mathbf{e}_{0}$ is $-y_{0}\mathbf{e}_{0}-y_{1}\mathbf{%
e}_{1}=-\mathbf{y}$; note that this is independent of our choice of $\mathbf{%
e}_{1}$. Thus we define 
\begin{equation}
-\mathbf{e}_{0}\odot \mathbf{y}=-\mathbf{y}  \label{eq:-e0_def}
\end{equation}
for all $\mathbf{y}\in \mathbb{S}$.

The definition (\ref{eq:-e0_def}) has a geometric interpretation analogous
to (\ref{eq:op_tmp1}). While (\ref{eq:op_tmp1}) defines $\mathbf{x}\odot 
\mathbf{y}$ as being the reflection of $\mathbf{e}_{0}\ast \mathbf{y}=J%
\mathbf{y}$ about the point $\mathbf{x}^{1/2}$, (\ref{eq:-e0_def}) defines $-%
\mathbf{e}_{0}\odot \mathbf{y}$ as being the reflection of $J\mathbf{y}$
across the subspace $\mathbf{V}$: 
\begin{equation*}
-\mathbf{e}_{0}\odot \mathbf{y}=J\mathbf{y}-2y_{0}\mathbf{e}_{0}.
\end{equation*}

The following definition summarizes the preceding discussion.

\begin{definition}
\label{defn:operation}Let $\mathbb{S}$ denote the unit sphere in a real
Hilbert space. Define a binary operation $\odot :\mathbb{S}\times \mathbb{S}%
\rightarrow \mathbb{S}$ by 
\begin{equation}
\mathbf{x}\odot \mathbf{y}=\left\{ 
\begin{array}{cc}
\left( 2P_{\mathbf{x}^{1/2}}-I\right) J\mathbf{y}, & \text{if }\mathbf{x}%
\neq -\mathbf{e}_{0} \\ 
-\mathbf{y}, & \text{if }\mathbf{x}=-\mathbf{e}_{0}
\end{array}
\right.  \label{eq:op}
\end{equation}
for $\mathbf{x},\mathbf{y}\in \mathbb{S}$.
\end{definition}

Although most of the properties of the magma $(\mathbb{S},\odot )$ will
follow from (\ref{eq:op}), there is another form for the operation $\odot $
which is also occasionally useful. For $\mathbf{x},\mathbf{y}\in \mathbb{S}$%
, $\mathbf{x}\neq -\mathbf{e}_{0}$, it easy to show by direct computation
that 
\begin{equation}
\mathbf{x}\odot \mathbf{y}=\langle \mathbf{x},J\mathbf{y}\rangle \mathbf{e}%
_{0}+\frac{y_{0}+\langle \mathbf{x},J\mathbf{y}\rangle }{1+x_{0}}\mathbf{x}%
_{\perp }+\mathbf{y}_{\perp }.  \label{eq:new_op}
\end{equation}

Fix a plane $\Pi \subset \mathbf{H}$ through $\mathbf{e}_{0}$ and the
origin, fix $\mathbf{e}_{1}\in \mathbb{S}\cap \Pi \cap \mathbf{V}$ and
consider the mapping $\mathbb{S}\cap \Pi \rightarrow \mathbb{C}:\mathbf{x}%
\longmapsto x_{0}+ix_{1}$, where $x_{1}=\mathbf{x}_{\perp }^{T}\mathbf{e}%
_{1} $. This mapping is injective, and its image is the unit circle $%
S^{1}\subset \mathbb{C}$. From (\ref{eq:circle}) we see that this mapping is
an isomorphism from the submagma $(\mathbb{S}\cap \Pi ,\odot )$ to the
circle group $(S^{1},\cdot )$ where the operation $\cdot $ refers to
multiplication of complex numbers. (This isomorphism depends on the choice
of $\mathbf{e}_{1}$, which fixes an orientation of the circle $\mathbb{S}%
\cap \Pi $.) We thus have the following result.

\begin{theorem}
\label{thm:subgroups}Let $\Pi \subset \mathbf{H}$ be a plane through $%
\mathbf{e}_{0}$ and the origin. Then $(\mathbb{S}\cap \Pi ,\odot )$ is an
abelian subgroup of $(\mathbb{S},\odot )$.
\end{theorem}

If we consider the case where $\mathbb{S}=S^{1}$, we have an immediately
corollary.

\begin{corollary}
\label{coro:circle_group}$(S^{1},\odot )=(S^{1},\cdot )$.
\end{corollary}

For $\mathbf{x}\in \mathbb{S}$, define 
\begin{equation}
\mathbf{x}^{-1}=J\mathbf{x}=x_{0}\mathbf{e}_{0}-\mathbf{x}_{\perp }\in 
\mathbb{S}.  \label{eq:inverse}
\end{equation}
In particular, $\mathbf{e}_{0}^{-1}=\mathbf{e}_{0}$ and $(-\mathbf{e}%
_{0})^{-1}=-\mathbf{e}_{0}$. The following is another corollary of Theorem 
\ref{thm:subgroups}.

\begin{corollary}
\label{coro:ident_inverses}

\begin{enumerate}
\item  The element $\mathbf{e}_{0}\in \mathbb{S}$ is a two-sided identity
for $(\mathbb{S},\odot )$.

\item  For each $\mathbf{x}\in \mathbb{S}$, $\mathbf{x}^{-1}$ is a two-sided
inverse of $\mathbf{x}$.

\item  $(\mathbb{S},\odot )$ is power-associative.
\end{enumerate}
\end{corollary}

Theorem \ref{thm:subgroups} and the definitions (\ref{eq:root}) and (\ref
{eq:inverse}) suggest another useful operation, namely an action of $\mathbb{%
R}$ on $\mathbb{S}$.

\begin{definition}
\label{defn:xt}Let $\mathbb{S}$ denote the unit sphere in a real Hilbert
space. For $\mathbf{x}\in \mathbb{S}$, $\mathbf{x}\neq -\mathbf{e}_{0}$, and
for $t\in \mathbb{R}$, define 
\begin{equation}
\mathbf{x}^{t}=\cos \left( t\cos ^{-1}x_{0}\right) \mathbf{e}_{0}+\sin
\left( t\cos ^{-1}x_{0}\right) \frac{\mathbf{x}_{\perp }}{\left\| \mathbf{x}%
_{\perp }\right\| }.  \label{eq:xt}
\end{equation}
(Here $\cos ^{-1}u=\arccos u$.) Also define $\mathbf{e}_{0}^{t}=\mathbf{e}%
_{0}$ for all $t\in \mathbb{R}$. Let $\omega :\mathbb{S}\backslash \{-%
\mathbf{e}_{0}\}\times \mathbb{R}\rightarrow \mathbb{S}$ denote the
operation $(\mathbf{x},t)\longmapsto \mathbf{x}^{t}$.
\end{definition}

Some properties easily follow from this definition. Once we have shown that
the magma $(\mathbb{S},\odot )$ is a left loop, it will follow from the next
result that $(\mathbb{S},\odot ,\omega )$ is a partial real odule (cf. \cite
{mik-sab}, XII.1.16).

\begin{theorem}
\label{thm:odule}Fix $\mathbf{x}\in \mathbb{S}$ with $\mathbf{x}\neq -%
\mathbf{e}_{0}$. The following hold.

\begin{enumerate}
\item  $\mathbf{x}^{1}=\mathbf{x}$;

\item  $(\mathbf{x}^{s})^{t}=\mathbf{x}^{st}$ for $0\leq s<\pi /\cos
^{-1}x_{0}$ and for all $t\in \mathbb{R}$;

\item  $\mathbf{x}^{t}\odot \mathbf{x}^{s}=\mathbf{x}^{t+s}$ for all $s,t\in 
\mathbb{R}$ such that $\mathbf{x}^{t},\mathbf{x}^{s},\mathbf{x}^{t+s}\neq -%
\mathbf{e}_{0}$.
\end{enumerate}
\end{theorem}

\begin{proof}
(1) is immediate from (\ref{eq:xt}) since $\sin (\cos ^{-1}x_{0})=\left\| 
\mathbf{x}_{\perp }\right\| $. To simplify further calculations, set $a=\cos
^{-1}x_{0}$. For (2), the restriction on $s$ implies $\cos ^{-1}(\cos
(sa))=sa$, and the result follows immediately. For (3), we compute 
\begin{equation*}
\langle \mathbf{x}^{s},J\mathbf{x}^{t}\rangle =\cos (sa)\cos (ta)-\sin
(ta)\sin (sa)=\cos ((s+t)a)
\end{equation*}
and 
\begin{eqnarray*}
\frac{\cos (sa)+\langle \mathbf{x}^{s},J\mathbf{x}^{t}\rangle }{1+\cos (ta)}%
\sin (ta) &=&\cos (sa)\sin (ta)-\frac{\sin ^{2}(ta)\sin (sa)}{1+\cos (ta)} \\
&=&\cos (sa)\sin (ta)-\sin (sa)(1-\cos (ta)) \\
&=&\sin ((t+s)a)-\sin (sa).
\end{eqnarray*}
Using these calculations in (\ref{eq:new_op}) gives the desired result.
\end{proof}

The restriction on $s$ in Theorem \ref{thm:odule}(2) is generically
unavoidable. For instance, if $\pi /\cos ^{-1}x_{0}<s<2\pi /\cos ^{-1}x_{0}$%
, then 
\begin{equation*}
(\mathbf{x}^{s})^{t}=\cos (t(2\pi -s\cos ^{-1}x_{0}))\mathbf{e}_{0}+\sin
\left( t(2\pi -s\cos ^{-1}x_{0})\right) \frac{\mathbf{x}_{\perp }}{\left\| 
\mathbf{x}_{\perp }\right\| }\neq \mathbf{x}^{s}{}^{t},
\end{equation*}
unless $t$ is an integer.

Some consequences of Definition \ref{defn:xt} and Theorem \ref{thm:odule}
are worth noting explicitly. First, $\mathbf{x}^{0}=\mathbf{e}_{0}$, and in
fact, for any integer $n\neq 0$, $\mathbf{x}^{n}=\mathbf{x}\odot \cdots
\odot \mathbf{x}$ (n factors). In particular, (\ref{eq:xt}) agrees with (\ref
{eq:inverse}) for $t=-1$. Also (\ref{eq:xt}) agrees with (\ref{eq:root}) for 
$t=1/2$.

We will consider additional properties of the magma $(\mathbb{S},\odot )$ in
the next section. For now we make a couple of remarks which anticipate later
results. We have seen that $(S^{1},\odot )$ is exactly the circle group.
However, $(S^{3},\odot )$ is not (isomorphic to) the group of unit
quaternions, and $(S^{7},\odot )$ is not (isomorphic to) the Moufang loop of
unit octonions. This is because $(\mathbb{S},\odot )$ turns out to be a left
loop, but not a loop for $\dim \mathbf{H}>2$. Also $(S^{15},\odot )$ is not
(isomorphic to) the left loop studied by Smith \cite{smith} using sedenian
multiplication; the latter left loop is not power-associative.

Finally, recall that the $2$-sphere $S^{2}$ may be identified with the \emph{%
Riemann sphere} $\mathbb{\hat{C}}=\mathbb{C}\cup \{\infty \}$ by
stereographic projection, where $\mathbf{e}_{0}$, the ``north pole'', is
mapped to $0\in \mathbb{C}$ and $-\mathbf{e}_{0}$, the ``south pole'', is
mapped to $\infty $. Using this mapping, it is straightforward to transfer
the operation $\odot $ to $\mathbb{\hat{C}}$. For $x,y\in \mathbb{\hat{C}}$,
we find that 
\begin{equation}
x\odot y=\left\{ 
\begin{array}{cc}
\dfrac{x+y}{1-\bar{x}y} & \text{if }x,y\neq \infty \\ 
-1/\bar{y} & \text{if }x=\infty \\ 
-1/\bar{x} & \text{if }y=\infty
\end{array}
\right. .  \label{eq:S2_op}
\end{equation}
Here $\bar{x}$ denotes the complex conjugate of $x$, and we are using the
usual conventions of complex arithmetic: $1/0=\infty $ and $1/\infty =0$.
For $x,y\neq \infty $, (\ref{eq:S2_op}) agrees with the local geodesic
operation on $\mathbb{\hat{C}}$ found by Nesterov and Sabinin \cite{nes-sab}%
. For the most part, results about $(\mathbb{\hat{C}},\odot )$ are just
special cases of results about $(\mathbb{S},\odot )$.

\section{Algebraic Structure}

It is clear from (\ref{eq:op}) that the left translation mappings $L_{%
\mathbf{x}}:\mathbb{S}\rightarrow \mathbb{S}$ defined by $L_{\mathbf{x}}%
\mathbf{y}=\mathbf{x}\odot \mathbf{y}$ are given by the restriction to $%
\mathbb{S}$ of bounded linear operators defined on all of $\mathbf{H}$.
Abusing notation slightly, we will denote these operators by $L_{\mathbf{x}}$
also. We have 
\begin{equation}
L_{\mathbf{x}}=\left\{ 
\begin{array}{cc}
\left( 2P_{\mathbf{x}^{1/2}}-I\right) J, & \text{if }\mathbf{x}\neq -\mathbf{%
e}_{0} \\ 
-I & \text{if }\mathbf{x}=-\mathbf{e}_{0}
\end{array}
\right. .  \label{eq:left_trans}
\end{equation}

\begin{theorem}
\label{thm:translation_props}Let $\mathbf{x}\in \mathbb{S}$ be given. The
following hold. 
\begin{eqnarray}
L_{\mathbf{x}} &\in &O(\mathbf{H}).  \label{eq:L_orthogonal} \\
JL_{\mathbf{x}} &=&L_{J\mathbf{x}}J.  \label{eq:JLx=LJxJ} \\
L_{\mathbf{x}}^{T} &=&L_{\mathbf{x}^{-1}}.  \label{eq:Ltrans=Lx'}
\end{eqnarray}
\end{theorem}

\begin{proof}
All three assertions are trivial for $\mathbf{x}=-\mathbf{e}_{0}$, and thus
we may assume $\mathbf{x}\neq -\mathbf{e}_{0}$. Eq. (\ref{eq:L_orthogonal})
follows from the computation 
\begin{equation*}
L_{\mathbf{x}}L_{\mathbf{x}}^{T}=4P_{\mathbf{x}^{1/2}}^{2}-4P_{\mathbf{x}%
^{1/2}}+I=I
\end{equation*}
using (\ref{eq:JJ=I}) and the fact that $P_{\mathbf{x}^{1/2}}$ is an
orthogonal projection. Now 
\begin{equation*}
J\mathbf{x}^{1/2}=\mathbf{x}^{-1/2}=(J\mathbf{x})^{1/2},
\end{equation*}
and thus 
\begin{equation*}
JL_{\mathbf{x}}=(2JP_{\mathbf{x}^{1/2}}-J)J=(2P_{(J\mathbf{x}%
)^{1/2}}-I)JJ=L_{J\mathbf{x}}J
\end{equation*}
using (\ref{eq:APa=PAaA}) (with $A=J$). This establishes (\ref{eq:JLx=LJxJ}%
). Finally, (\ref{eq:Ltrans=Lx'}) follows from (\ref{eq:JLx=LJxJ}), (\ref
{eq:JJ=I}), and the symmetry of $P_{\mathbf{x}^{1/2}}$: 
\begin{equation*}
L_{\mathbf{x}}^{T}=JL_{\mathbf{x}}J=L_{\mathbf{x}^{-1}}.
\end{equation*}
This completes the proof.
\end{proof}

Let $\mathrm{Aut}(\mathbb{S},\odot )$ denote the group of continuous
automorphisms of $(\mathbb{S},\odot )$.

\begin{theorem}
\label{thm:automorphisms}$\mathrm{Aut}(\mathbb{S},\odot )=O(\mathbf{V})$.
\end{theorem}

We omit the (long) proof, but note that the containment $O(\mathbf{V}%
)\subseteq \mathrm{Aut}(\mathbb{S},\odot )$ is clear from (\ref{eq:new_op}).

We will identify each left inner mapping $L(\mathbf{x},\mathbf{y})=L_{%
\mathbf{x}\odot \mathbf{y}}^{-1}L_{\mathbf{x}}L_{\mathbf{y}}$ on $\mathbb{S}$
with its extension to a bounded linear operator on $\mathbf{H}$.

\begin{theorem}
\label{thm:left_inner_props}For all $\mathbf{x},\mathbf{y}\in \mathbb{S}$, $%
L(\mathbf{x},\mathbf{y})\in O(\mathbf{V})$.
\end{theorem}

\begin{proof}
Since $O(\mathbf{H})$ is a group, $L(\mathbf{x},\mathbf{y})=L_{\mathbf{x}%
\odot \mathbf{y}}^{-1}L_{\mathbf{x}}L_{\mathbf{y}}\in O(\mathbf{H})$ by (\ref
{eq:L_orthogonal}). We have $L(\mathbf{x},\mathbf{y})\mathbf{e}_{0}=L_{%
\mathbf{x}\odot \mathbf{y}}^{-1}\left( \mathbf{x}\odot \mathbf{y}\right) =%
\mathbf{e}_{0}$. In view of (\ref{eq:O(V)}), this establishes the result.
\end{proof}

Putting the preceding discussion together, we have identified the class of
left loops to which $(\mathbb{S},\odot )$ belongs.

\begin{corollary}
\label{coro:kikkawa}$(\mathbb{S},\odot )$ is a Kikkawa left loop.
\end{corollary}

\begin{proof}
The left inverse property follows from (\ref{eq:L_orthogonal}) and (\ref
{eq:Ltrans=Lx'}). The automorphic inverse property follows from (\ref
{eq:JLx=LJxJ}). The $A_{l}$ property follows from Theorems \ref
{thm:left_inner_props} and \ref{thm:automorphisms}.
\end{proof}

Not only do operators in $O(\mathbf{V})$ preserve the operation $\odot $,
they also preserve the scalar multiplication $\omega $ given by (\ref{eq:xt}%
). Indeed, for all $\mathbf{x}\in \mathbb{S}$, $\mathbf{x}\neq -\mathbf{e}%
_{0}$, for all $t\in \mathbb{R}$, and for all $A\in O(\mathbf{V})$, we have 
\begin{equation}
(A\mathbf{x})^{t}=A\mathbf{x}^{t}  \label{eq:(Ax)t=Axt}
\end{equation}
as an obvious consequence of the definitions. From (\ref{eq:(Ax)t=Axt}) and
Theorem \ref{thm:left_inner_props}, we have 
\begin{equation}
(L(\mathbf{x},\mathbf{y})\mathbf{z})^{t}=L(\mathbf{x},\mathbf{y})\mathbf{z}%
^{t}  \label{eq:2ndAl}
\end{equation}
for all $\mathbf{x},\mathbf{y},\mathbf{z}\in \mathbb{S}$, $\mathbf{z}\neq -%
\mathbf{e}_{0}$, $t\in \mathbb{R}$. Thus the $\mathbb{R}$-odule $(\mathbb{S}%
,\odot ,\omega )$ satisfies Sabinin's so-called ``second $A_{l}$ property''
(cf. \cite{mik-sab}, XII.3.18).

Kiechle has shown that there exist Kikkawa left loops which are not loops ( 
\cite{kiechle}, 11.3.2). Here we will show that $(\mathbb{S},\odot )$ is not
a loop. For $\mathbf{a}\in \mathbb{S}$, consider the equation 
\begin{equation}
\mathbf{x}\odot \mathbf{a}=-\mathbf{a}^{-1}  \label{eq:x.a=-a'}
\end{equation}
in $\mathbb{S}$. From Theorem \ref{thm:subgroups}, $\mathbf{x}=-\mathbf{a}%
^{-2}$ is a solution. If $\mathbf{a}=\pm \mathbf{e}_{0}$, then clearly $%
\mathbf{x}=-\mathbf{a}^{-2}=-\mathbf{e}_{0}$ is the only solution.
Conversely, if $\mathbf{x}=-\mathbf{e}_{0}$ is a solution, then $-\mathbf{a}%
=-\mathbf{a}^{-1}$, i.e., $\mathbf{a}=\pm \mathbf{e}_{0}$. In general, (\ref
{eq:x.a=-a'}) has solutions other than $\mathbf{x}=-\mathbf{a}^{-2}$.

\begin{theorem}
\label{thm:not_loop}For all $\mathbf{a}\in \mathbb{S}$, $\mathbf{a}\neq \pm 
\mathbf{e}_{0}$, the equation (\ref{eq:x.a=-a'}) has the solution set 
\begin{equation}
\mathbf{X}=\left[ -\mathbf{a}^{-2}+(\mathbf{a}^{-1})^{\perp }\right] \cap 
\mathbb{S}\backslash \left\{ -\mathbf{e}_{0}\right\} .  \label{eq:X}
\end{equation}
If $\dim \mathbf{H}>2$, then $\mathbf{X}$ is not a singleton.
\end{theorem}

\begin{proof}
By the preceding discussion, $\mathbf{x}\neq -\mathbf{e}_{0}$. Thus (\ref
{eq:x.a=-a'}) becomes $(2P_{\mathbf{x}^{1/2}}-I)J\mathbf{a}=-\mathbf{a}^{-1}$%
, i.e., $P_{\mathbf{x}^{1/2}}\mathbf{a}^{-1}=\mathbf{0}$. This is equivalent
to $\langle \mathbf{x}^{1/2},\mathbf{a}^{-1}\rangle =0$, and thus to $%
\langle \mathbf{e}_{0}+\mathbf{x},\mathbf{a}^{-1}\rangle =0$, or $\langle 
\mathbf{x},\mathbf{a}^{-1}\rangle =-a_{0}$, with the constraint $\mathbf{x}%
\in \mathbb{S}\backslash \left\{ -\mathbf{e}_{0}\right\} $. An obvious
particular solution is $\mathbf{x}=-\mathbf{a}^{-2}$. Thus the solution set
in $\mathbf{H}$ is the affine subspace $-\mathbf{a}^{-2}+(\mathbf{a}%
^{-1})^{\perp }$, and hence the solution set in $\mathbb{S}$ is $\mathbf{X}$%
, as claimed. The nonuniqueness assertion follows from geometric
considerations: $\mathbf{X}$ is a punctured sphere of codimension $2$ in $%
\mathbf{H}$ with center $-a_{0}\mathbf{a}^{-1}$ and containing $-\mathbf{a}%
^{-2}$. The antipode of $-\mathbf{a}^{-2}$ is the removed point $-\mathbf{e}%
_{0}$. If $\mathbf{X}$ is a singleton, then either $(\mathbf{a}^{-1})^{\perp
}\cap \mathbb{S}=\left\{ -\mathbf{a}^{-2},-\mathbf{e}_{0}\right\} $, in
which case $\dim \mathbf{H}\leq 2$, or $(\mathbf{a}^{-1})^{\perp }$ is
tangent to $\mathbb{S}$ at $-\mathbf{a}^{-2}$, in which case $-\mathbf{a}%
^{-2}=\pm \mathbf{a}^{-1}$ or $\mathbf{a}=\pm \mathbf{e}_{0}$, a
contradiction.
\end{proof}

We now turn to properties of $(\mathbb{S},\odot )$ which fail to hold
everywhere, but which hold ``almost everywhere''. First we consider the left
power alternative property: $\mathbf{x}^{s}\odot (\mathbf{x}^{t}\odot 
\mathbf{y})=\mathbf{x}^{s+t}\odot \mathbf{y}$. Special cases are the left
inverse property ($s=\pm 1$, $t=\mp 1$) (which we already know to hold
everywhere) and the left alternative property ($s=t=1$).

\begin{theorem}
\label{thm:lpa1}For all $\mathbf{x}\in \mathbb{S}$ with $\mathbf{x}\neq -%
\mathbf{e}_{0}$ and for all $s,t\in \mathbb{R}$ such that $\mathbf{x}^{s},%
\mathbf{x}^{t}\neq -\mathbf{e}_{0}$, 
\begin{equation}
L_{\mathbf{x}^{s}}L_{\mathbf{x}^{t}}=(2P_{\mathbf{x}^{(s+t)/2}}-I)J.
\label{eq:lpa1}
\end{equation}
In addition, $\mathbf{x}^{s+t}\neq -\mathbf{e}_{0}$ if and only if 
\begin{equation}
L_{\mathbf{x}^{s}}L_{\mathbf{x}^{t}}=L_{\mathbf{x}^{s+t}}.  \label{eq:lpa2}
\end{equation}
\end{theorem}

\begin{proof}
(Sketch) Fix $\mathbf{x}\in \mathbb{S}$ with $\mathbf{x}\neq -\mathbf{e}_{0}$%
. To simplify notation, set $C(t)=\cos ((t/2)\cos ^{-1}x_{0})$ and $%
S(t)=\sin ((t/2)\cos ^{-1}x_{0})$. For $s,t\in \mathbb{R}$ such that $%
\mathbf{x}^{s},\mathbf{x}^{t}\neq -\mathbf{e}_{0}$, we have 
\begin{equation}
L_{\mathbf{x}^{s}}L_{\mathbf{x}^{t}}J=4P_{\mathbf{x}^{t/2}}P_{\mathbf{x}%
^{s/2}}-2(JP_{\mathbf{x}^{s/2}}+P_{\mathbf{x}^{t/2}}J)+J  \label{eq:ugh}
\end{equation}
using (\ref{eq:left_trans}). Now 
\begin{eqnarray*}
P_{\mathbf{x}^{t/2}}P_{\mathbf{x}^{s/2}} &=&\langle \mathbf{x}^{t/2},J%
\mathbf{x}^{s/2}\rangle \mathbf{x}^{t/2}(\mathbf{x}^{s/2})^{T} \\
&=&(C(t)C(s)-S(t)S(s))\mathbf{x}^{t/2}(\mathbf{x}^{s/2})^{T} \\
&=&C(t+s)\mathbf{x}^{t/2}(\mathbf{x}^{s/2})^{T}.
\end{eqnarray*}
Using (\ref{eq:I+J=2e0e0}) and (\ref{eq:xt}), we may expand (\ref{eq:ugh})
into an expression of the form 
\begin{equation*}
L_{\mathbf{x}^{s}}L_{\mathbf{x}^{t}}J=\alpha \mathbf{e}_{0}\mathbf{e}%
_{0}^{T}+\frac{\beta }{\left\| \mathbf{x}_{\perp }\right\| }\mathbf{e}_{0}%
\mathbf{x}_{\perp }^{T}+\frac{\gamma }{\left\| \mathbf{x}_{\perp }\right\| }%
\mathbf{x}_{\perp }\mathbf{e}_{0}^{T}+\frac{\delta }{\left\| \mathbf{x}%
_{\perp }\right\| ^{2}}\mathbf{x}_{\perp }\mathbf{x}_{\perp }^{T}-I,
\end{equation*}
where the coefficients $\alpha ,\beta ,\gamma ,\delta $ are to be
determined. We find, for instance, 
\begin{eqnarray*}
\alpha &=&4C(t+s)C(t)C(s)-2C(s)^{2}-2S(t)^{2}+2 \\
&=&2C(t+s)(C(t+s)+C(t-s))-2(C(s)^{2}C(t)^{2}-S(s)^{2}S(t)^{2}) \\
&=&2C(t+s)^{2}.
\end{eqnarray*}
Similar calculations give 
\begin{eqnarray*}
\beta &=&\gamma =2C(t+s)S(t+s) \\
\delta &=&2S(t+s)^{2}.
\end{eqnarray*}
Thus $L_{\mathbf{x}^{s}}L_{\mathbf{x}^{t}}J=$ 
\begin{equation*}
2\left( C(t+s)\mathbf{e}_{0}+S(t+s)\frac{\mathbf{x}_{\perp }}{\left\| 
\mathbf{x}_{\perp }\right\| }\right) \left( C(t+s)\mathbf{e}_{0}+S(t+s)\frac{%
\mathbf{x}_{\perp }}{\left\| \mathbf{x}_{\perp }\right\| }\right) ^{T}-I.
\end{equation*}
This and (\ref{eq:xt}) give the desired result. The remaining assertion is
clear.
\end{proof}

Next we consider the cases ruled out in the hypotheses of Theorem \ref
{thm:lpa1}.

\begin{theorem}
\label{thm:lpa2}If $\dim \mathbf{H}>2$, then for all $\mathbf{x}\in \mathbb{S%
}$ with $\mathbf{x}\neq \pm \mathbf{e}_{0}$, $L_{-e_{0}}L_{\mathbf{x}}=L_{%
\mathbf{x}}L_{-\mathbf{e}_{0}}\neq L_{-\mathbf{x}}.$
\end{theorem}

\begin{proof}
Fix $\mathbf{x}\in \mathbb{S}$ with $\mathbf{x}\neq \pm \mathbf{e}_{0}$, and
assume $L_{\mathbf{x}}L_{-\mathbf{e}_{0}}=L_{-\mathbf{x}}$. Then $L_{\mathbf{%
x}^{-1}}L_{-\mathbf{x}}=-I$. Choose $\mathbf{y}\in \mathbb{S}\cap \mathbf{V}%
\cap \mathbf{x}^{\perp }$. Then $\mathbf{y}\in (J\mathbf{x})^{\perp }$.
Suppressing details, we have $-\mathbf{y}=L_{\mathbf{x}^{-1}}L_{-\mathbf{x}}%
\mathbf{y}=-L_{\mathbf{x}^{-1}}J\mathbf{y}=\mathbf{y}$, a contradiction.
\end{proof}

The remaining case is clear: 
\begin{equation}
L_{-\mathbf{e}_{0}}L_{-\mathbf{e}_{0}}=I=L_{\mathbf{e}_{0}}=L_{(-\mathbf{e}%
_{0})^{2}}.  \label{eq:lap-easy}
\end{equation}
Putting together Theorems \ref{thm:lpa1} and \ref{thm:lpa2} with (\ref
{eq:lap-easy}), we easily characterize those points in $\mathbb{S}$ for
which the left alternative property holds.

\begin{corollary}
\label{coro:alternative}Assume $\dim \mathbf{H}>2$. For all $\mathbf{x}\in 
\mathbb{S}$, $L_{\mathbf{x}}^{2}=L_{\mathbf{x}^{2}}$ if and only if $\mathbf{%
x}\notin \mathbf{V}$.
\end{corollary}

In particular, the left alternative property holds ``almost everywhere'',
since it fails only on the set $\mathbb{S}\cap \mathbf{V}$ of ``measure
zero''. In the finite dimensional case, we may remove the quotation marks.

A globally smooth Kikkawa left loop necessarily satisfies the left
alternative property (\cite{kikkawa}, Lemma 6.2. The result is stated for
loops, but the proof clearly works for left loops.) It follows from
Corollary \ref{coro:alternative} that the multiplication $\odot $ on $%
\mathbb{S}$ is not globally smooth. We will examine the continuity of $\odot 
$ more closely in the next section.

Next we consider the Bol identity $L_{\mathbf{x}}L_{\mathbf{y}}L_{\mathbf{x}%
}=L_{\mathbf{x}\odot (\mathbf{y}\odot \mathbf{x})}$. Since the left
alternative law does not hold everywhere, the Bol identity cannot hold
everywhere either. Observe that the identity is trivial for $\mathbf{x}=\pm 
\mathbf{e}_{0}$.

\begin{theorem}
\label{thm:no_Bol}Assume $\dim \mathbf{H}>2$. Let $\mathbf{x}\in \mathbb{S}$%
, $\mathbf{x}\neq \pm \mathbf{e}_{0}$, be given. If $\mathbf{y}=-\mathbf{e}%
_{0}$ or if $\mathbf{y}\neq -\mathbf{e}_{0},-\mathbf{x}^{-2}$ and $\mathbf{y}%
\odot \mathbf{x}=-\mathbf{x}^{-1}$, then $L_{\mathbf{x}}L_{\mathbf{y}}L_{%
\mathbf{x}}\neq L_{\mathbf{x}\odot (\mathbf{y}\odot \mathbf{x})}$.
\end{theorem}

\begin{proof}
We have 
\begin{equation*}
L_{\mathbf{x}}L_{-\mathbf{e}_{0}}L_{\mathbf{x}}=-L_{\mathbf{x}}^{2}=-(2P_{%
\mathbf{x}}-I)J
\end{equation*}
by (\ref{eq:lpa1}), and 
\begin{equation*}
L_{\mathbf{x}\odot (-\mathbf{e}_{0}\odot \mathbf{x})}=L_{\mathbf{x}\odot (-%
\mathbf{x})}=L_{-\mathbf{x}^{2}}=(2P_{(-\mathbf{x}^{2})^{1/2}}-I)J.
\end{equation*}
Now choose $\mathbf{z}\in \mathbf{V}\cap \mathbf{x}^{\perp }$. We have $L_{%
\mathbf{x}}L_{-\mathbf{e}_{0}}L_{\mathbf{x}}\mathbf{z}=J\mathbf{z}=-\mathbf{z%
}$ and $L_{\mathbf{x}\odot (-\mathbf{e}_{0}\odot \mathbf{x})}\mathbf{z}=-J%
\mathbf{z}=\mathbf{z}$. Thus $L_{\mathbf{x}}L_{-\mathbf{e}_{0}}L_{\mathbf{x}%
}=L_{\mathbf{x}\odot (-\mathbf{e}_{0}\odot \mathbf{x})}$ if and only if $%
\mathbf{V}\cap \mathbf{x}^{\perp }=\left\{ \mathbf{0}\right\} $, i.e., if
and only if $\dim \mathbf{H}\leq 2$.

Now assume $\mathbf{y}\neq -\mathbf{e}_{0}$ and $\mathbf{y}\odot \mathbf{x}=-%
\mathbf{x}^{-1}$ (i.e., $\mathbf{y}\in \mathbf{X}$ is given by (\ref{eq:X}%
)). Then $\mathbf{x}\odot (\mathbf{y}\odot \mathbf{x})=-\mathbf{e}_{0}$, and
thus $L_{\mathbf{x}}^{-1}L_{\mathbf{x}\odot (\mathbf{y}\odot \mathbf{x})}L_{%
\mathbf{x}}^{-1}=-L_{\mathbf{x}}^{-2}$. If the Bol identity holds, then $L_{%
\mathbf{y}}=-L_{\mathbf{x}}^{-2}$. Applying both sides to $\mathbf{e}_{0}$,
we have $\mathbf{y}=-\mathbf{x}^{-2}$.
\end{proof}

The conditions described in the previous theorem turn out to be the only
obstructions to the Bol identity.

\begin{theorem}
\label{thm:bol}Let $\mathbf{x},\mathbf{y}\in \mathbb{S}$, $\mathbf{x}\neq
\pm \mathbf{e}_{0}$, $\mathbf{y}\neq -\mathbf{e}_{0}$ and $\mathbf{y}\odot 
\mathbf{x}\neq -\mathbf{x}^{-1}$. Then $L_{\mathbf{x}}L_{\mathbf{y}}L_{%
\mathbf{x}}=L_{\mathbf{x}\odot (\mathbf{y}\odot \mathbf{x})}$.
\end{theorem}

We omit the proof, but outline the idea. Applying both sides of the Bol
identity to $\mathbf{e}_{0}$ gives an obvious equality, so it is only
necessary to check that both sides agree when applied to an arbitrary
element of $\mathbf{V}$. The hypotheses guarantee that all the left
translations in the identity have the form $(2P_{\mathbf{u}^{1/2}}-I)J$, so
all that remains is a calculation. This is tedious, but straightforward.

\section{Continuity and Smoothness}

We have already noted in the previous section that the operation $\odot $
cannot be globally smooth on $\mathbb{S}\times \mathbb{S}$ because $(\mathbb{%
S},\odot )$ does not everywhere satisfy the left alternative property.

For $\mathbf{x}\in \mathbb{S}$ with $\mathbf{x}\neq -\mathbf{e}_{0}$, we
have 
\begin{equation}
\lim_{t\rightarrow \pi /\cos ^{-1}x_{0}}\mathbf{x}^{t}=-\mathbf{e}_{0}
\label{eq:limxt=-e0}
\end{equation}
On the other hand, $\lim_{t\rightarrow \pi /\cos ^{-1}x_{0}}\mathbf{x}^{t/2}=%
\mathbf{x}_{\perp }/\left\| \mathbf{x}_{\perp }\right\| $, and thus 
\begin{equation}
\lim_{t\rightarrow \pi /\cos ^{-1}x_{0}}L_{\mathbf{x}^{t}}=\left( 2\frac{%
\mathbf{x}_{\perp }\mathbf{x}_{\perp }^{T}}{\left\| \mathbf{x}_{\perp
}\right\| ^{2}}-I\right) J.  \label{eq:limLxt}
\end{equation}
Now for $\mathbf{y}\in \mathbb{S}$, $\left( (2/\left\| \mathbf{x}_{\perp
}\right\| ^{2})\mathbf{x}_{\perp }\mathbf{x}_{\perp }^{t}-I\right) J\mathbf{y%
}=-\mathbf{y}$ if and only if $-\left\| \mathbf{x}_{\perp }\right\| ^{2}%
\mathbf{y}_{\perp }=(\mathbf{x}_{\perp }^{T}\mathbf{y}_{\perp })\mathbf{x}%
_{\perp }$ i.e., if and only $\mathbf{x}$, $\mathbf{y}$, and $\mathbf{e}_{0}$
are coplanar with $\mathbf{0}$. Therefore, for $\mathbf{x},\mathbf{y}\in 
\mathbb{S}$ with $\mathbf{x}\neq -\mathbf{e}_{0}$, and $\mathbf{x},\mathbf{y}%
,\mathbf{e}_{0}$ not coplanar with $\mathbf{0}$, 
\begin{equation}
\lim_{t\rightarrow \pi /\cos ^{-1}x_{0}}\mathbf{x}^{t}\odot \mathbf{y}\neq -%
\mathbf{e}_{0}\odot \mathbf{y}.  \label{eq:limbad}
\end{equation}
For $\mathbf{y}\in \mathbb{S}$, define the right translation operator $R_{%
\mathbf{y}}:\mathbb{S}\rightarrow \mathbb{S}:\mathbf{x}\longmapsto \mathbf{x}%
\odot \mathbf{y}$. We have established the following.

\begin{theorem}
\label{thm:discontinuity}For each $\mathbf{y}\in \mathbb{S}$, $R_{\mathbf{y}%
}:\mathbb{S}\rightarrow \mathbb{S}$ has a nonremovable discontinuity at $-%
\mathbf{e}_{0}$.
\end{theorem}

This turns out to be the only topological, or even analytic, deficiency.

\begin{proposition}
\label{prop:analytic}The following hold.

\begin{enumerate}
\item  For each $\mathbf{x}\in \mathbb{S}$, $L_{\mathbf{x}}:\mathbb{S}%
\rightarrow \mathbb{S}$ is analytic.

\item  For each $\mathbf{y}\in \mathbb{S}$, $R_{\mathbf{y}}:\mathbb{S}%
\backslash \{-\mathbf{e}_{0}\}\rightarrow \mathbb{S}$ is analytic.

\item  $\odot :\mathbb{S}\backslash \{-\mathbf{e}_{0}\}\times \mathbb{S}%
\rightarrow \mathbb{S}$ is analytic.
\end{enumerate}
\end{proposition}

\begin{proof}
(Sketch) Note that (1) is immediate since $L_{\mathbf{x}}$ is the
restriction of a bounded linear operator to $\mathbb{S}$. For (2), it is
clear from (\ref{eq:new_op}) that $R_{\mathbf{y}}$ ($\mathbf{y}\in \mathbb{S}
$) has an analytic extension to $\mathbf{H}\backslash \left\{ -\mathbf{e}%
_{0}\right\} $. Similarly, (\ref{eq:new_op}) shows that $\odot $ has an
analytic extension to $\mathbf{H}\backslash \left\{ -\mathbf{e}_{0}\right\}
\times \mathbf{H}$, and this implies (3).
\end{proof}

\section{Semidirect Product Structure}

In this section we examine the relationship between the left loop structure
of $(\mathbb{S},\odot )$ and the structure of the orthogonal group $O(%
\mathbf{H})$ relative to the subgroup $O(\mathbf{V})$. This discussion is a
particular case of the general theory of semidirect products of left loops
with groups. This theory was worked out in general by Sabinin \cite
{sabinin72} (see also \cite{mik-sab}, \cite{sab-book}), and was later
rediscovered in the particular case of $A_{l}$ left loops with LIP by
Kikkawa \cite{kikkawa} and Ungar \cite{ungar-wag}. A survey with recent
extensions can be found in \cite{kj}.

Since $(\mathbb{S},\odot )$ is an $A_{l}$ left loop, we may form its \emph{%
standard semidirect product} with $\mathrm{Aut}(\mathbb{S},\odot )$. This is
the group denoted by $\mathbb{S}\rtimes\mathrm{Aut}(\mathbb{S},\odot )$,
consisting of the set $\mathbb{S}\times \mathrm{Aut}(\mathbb{S},\odot )$
with multiplication defined by 
\begin{equation}
(\mathbf{x},A)(\mathbf{y},B)=(\mathbf{x}\odot A\mathbf{y},L(\mathbf{x},A%
\mathbf{y})AB).  \label{eq:semi-mult}
\end{equation}

Define 
\begin{equation}
L(\mathbb{S})=\left\{ L_{\mathbf{x}}:\mathbf{x}\in \mathbb{S}\right\}
\label{eq:L(S)}
\end{equation}
to be the set of all left translations. Then $L(\mathbb{S})\subset O(\mathbf{%
H})$ by (\ref{eq:L_orthogonal}), and $L(\mathbb{S})$ is a left transversal
of the subgroup $O(\mathbf{V})$. Indeed, for $A\in O(\mathbf{H})$, let $%
\mathbf{u}=A\mathbf{e}_{0}$. Then $U=L_{\mathbf{u}}^{-1}A\in O(\mathbf{V})$.
If $L_{\mathbf{u}}U=L_{\mathbf{v}}V$ with $\mathbf{v}\in \mathbb{S}$, $V\in
O(\mathbf{H})$, then applying both sides to $\mathbf{e}_{0}$, we obtain $%
\mathbf{u}=\mathbf{v}$, and thus $U=V$. Thus $A=L_{\mathbf{u}}U$ is a unique
factorization of $A$ into a left translation in $L(\mathbb{S})$ and an
operator in $O(\mathbf{V})$.

The transversal decomposition $O(\mathbf{H})=L(\mathbb{S})O(\mathbf{V})$
defines a left loop structure on $L(\mathbb{S})$ itself by projection: $L_{%
\mathbf{x}}\odot L_{\mathbf{y}}=L_{\mathbf{u}}$ where $L_{\mathbf{u}}$ is
the unique representative of $L_{\mathbf{x}}L_{\mathbf{y}}$ in $L(\mathbb{S}%
) $. Thus $O(\mathbf{H})$ is an \emph{internal semidirect product} of the
left loop $(L(\mathbb{S}),\odot )$ by the subgroup $O(\mathbf{V})$.

Now for $\mathbf{x},\mathbf{y}\in \mathbb{S}$, $A,B\in O(\mathbf{V})$, $L_{%
\mathbf{x}}AL_{\mathbf{y}}B=L_{\mathbf{x}\odot A\mathbf{y}}L(\mathbf{x},A%
\mathbf{y})AB$. Thus the left loops $(\mathbb{S},\odot )$ and $(L(\mathbb{S}%
),\odot )$ are just isomorphic copies (take $A=I$), and since $O(\mathbf{V}%
)\cong \mathrm{Aut}(\mathbb{S},\odot )$ (Theorem \ref{thm:automorphisms}),
we also have the following result.

\begin{theorem}
$O(\mathbf{H})\cong \mathbb{S}\rtimes O(\mathbf{V}).$
\end{theorem}

Interestingly, the transversal $L(\mathbb{S})$ in $O(\mathbf{H})$ is \emph{%
not} connected. The element $L_{-\mathbf{e}_{0}}=-I$ is an isolated point.
An intuitive argument can be seen as follows (cf. (\ref{eq:limbad})). For
any $\mathbf{u}\in \mathbb{S}$ with $\mathbf{u}\neq \pm \mathbf{e}_{0}$, 
\begin{equation*}
\lim_{\mathbf{x}\rightarrow \mathbf{u}}L_{\mathbf{x}}=\lim_{\mathbf{x}%
\rightarrow \mathbf{u}}(2P_{\mathbf{x}}-I)J=(2P_{\mathbf{u}}-I)J.
\end{equation*}
If $(2P_{\mathbf{u}}-I)J=-I$, then $2P_{\mathbf{u}}=I-J$. Apply this to any $%
\mathbf{v}\in \mathbb{S}\cap \mathbf{V}\cap \mathbf{u}^{\perp }$. Then $%
\mathbf{0}=(I-J)\mathbf{v}=2\mathbf{v}$, a contradiction. Thus $\lim_{%
\mathbf{x}\rightarrow \mathbf{u}}L_{\mathbf{x}}\neq L_{-\mathbf{e}_{0}}$.
This argument can easily be made more precise; we omit the details here. The
nonconnectedness of $L(\mathbb{S})$ shows once again that the mapping $L:%
\mathbb{S}\rightarrow O(\mathbf{H}):\mathbf{x}\longmapsto L_{\mathbf{x}}$ is
not everywhere continuous.

\section{Spherical Geometry}

We now make a few remarks about the relationship between the Kikkawa left
loop structure of $(\mathbb{S},\odot )$ and spherical geometry on $\mathbb{S}
$. This relationship is analogous to that between the $B$-loop structure of
the unit ball in a Hilbert space and hyperbolic geometry; see \cite
{ungar-axioms} and references therein.

On $\mathbb{S}$ we define a norm by 
\begin{equation}
\left\| \mathbf{x}\right\| _{s}=\cos ^{-1}x_{0}=\cos ^{-1}\langle \mathbf{x},%
\mathbf{e}_{0}\rangle  \label{eq:norm}
\end{equation}
for $\mathbf{x}\in \mathbb{S}$. Note the formal similarities between the
following properties of $\left\| \cdot \right\| _{s}:\mathbb{S}\rightarrow 
\mathbb{R}$, and properties of norms in vector spaces.

\begin{theorem}
\label{thm:norm_props}The following hold.

\begin{enumerate}
\item  For all $\mathbf{x}\in \mathbb{S}$, $0\leq \left\| \mathbf{x}\right\|
_{s}\leq \pi $. In addition, $\left\| \mathbf{x}\right\| _{s}=0$ if and only
if $\mathbf{x}=\mathbf{e}_{0}$, and $\left\| \mathbf{x}\right\| _{s}=\pi $
if and only if $\mathbf{x}=-\mathbf{e}_{0}$.

\item  For all $\mathbf{x}\in \mathbb{S}$, $\mathbf{x}\neq -\mathbf{e}_{0}$, 
$t\in \mathbb{R}$, $\left\| \mathbf{x}^{t}\right\| _{s}=\left| t\right|
\left\| \mathbf{x}\right\| _{s}\ \mathrm{mod}\pi $.

\item  (Triangle inequality) For all $\mathbf{x},\mathbf{y}\in \mathbb{S}$, 
\begin{equation}
\left\| \mathbf{x}\odot \mathbf{y}\right\| _{s}\leq \pi -\left| \left\| 
\mathbf{x}\right\| _{s}+\left\| \mathbf{y}\right\| _{s}-\pi \right| \leq
\left\| \mathbf{x}\right\| _{s}+\left\| \mathbf{y}\right\| _{s}.
\label{eq:tri}
\end{equation}
The first inequality is an equality if and only if $\mathbf{x}$, $\mathbf{y}$%
, and $\mathbf{e}_{0}$ are coplanar with $\mathbf{0}$.

\item  (Invariance under $O(\mathbf{V})$) For all $\mathbf{x}\in \mathbb{S}$%
, $A\in O(\mathbf{V})$, $\left\| A\mathbf{x}\right\| _{s}=\left\| \mathbf{x}%
\right\| _{s}$.
\end{enumerate}
\end{theorem}

\begin{proof}
(1), (2), and (4) are easy consequences of the definitions. For (3), we use (%
\ref{eq:new_op}) and the Cauchy-Schwarz inequality to compute 
\begin{eqnarray*}
\cos \left\| \mathbf{x}\odot \mathbf{y}\right\| _{s} &=&\langle \mathbf{x}%
\odot \mathbf{y},\mathbf{e}_{0}\rangle =\langle \mathbf{x},J\mathbf{y}%
\rangle =x_{0}y_{0}-\langle \mathbf{x}_{\perp },\mathbf{y}_{\perp }\rangle \\
&\geq &x_{0}y_{0}-\left\| \mathbf{x}_{\perp }\right\| \left\| \mathbf{y}%
_{\perp }\right\| \\
&=&\cos (\cos ^{-1}x_{0})\cos (\cos ^{-1}y_{0})-\sin \left( \cos
^{-1}x_{0}\right) \sin (\cos ^{-1}y_{0}) \\
&=&\cos (\left\| \mathbf{x}\right\| _{s}+\left\| \mathbf{y}\right\| _{s}).
\end{eqnarray*}
By (1), we have 
\begin{equation*}
\left\| \mathbf{x}\odot \mathbf{y}\right\| _{s}\leq \cos ^{-1}(\cos (\left\| 
\mathbf{x}\right\| _{s}+\left\| \mathbf{y}\right\| _{s}))=\pi -\left|
\left\| \mathbf{x}\right\| _{s}+\left\| \mathbf{y}\right\| _{s}-\pi \right| ,
\end{equation*}
since $0\leq \left\| \mathbf{x}\right\| _{s}+\left\| \mathbf{y}\right\|
_{s}\leq 2\pi $. If equality holds, then the Cauchy-Schwarz equality implies
that $\mathbf{x}_{\perp }$ and $\mathbf{y}_{\perp }$ are parallel, which is
equivalent to the desired result. The remaining inequality of (\ref{eq:tri})
is clear.
\end{proof}

Next we use the operation $\odot $ and the norm to define a distance
function 
\begin{equation}
d_{s}(\mathbf{x},\mathbf{y})=\left\| \mathbf{x}\odot \mathbf{y}^{-1}\right\|
_{s}  \label{eq:distance}
\end{equation}
for $\mathbf{x},\mathbf{y}\in \mathbb{S}$. This definition is analogous to
the usual relationship between norms and distance functions in vector
spaces. By (\ref{eq:new_op}), we have $\langle \mathbf{x}\odot \mathbf{y}%
^{-1},\mathbf{e}_{0}\rangle =\langle \mathbf{x},J\mathbf{y}^{-1}\rangle
=\langle \mathbf{x},\mathbf{y}\rangle $, and hence $d_{s}(\mathbf{x},\mathbf{%
y})=\cos ^{-1}\langle \mathbf{x},\mathbf{y}\rangle $. This shows that (\ref
{eq:distance}) is equivalent to the usual definition of the spherical
distance function.

Before establishing properties of the distance function, we require a lemma.
The following result was established by Ungar for Bruck loops \cite
{ungar-disk}. Here we extend it to Kikkawa left loops.

\begin{lemma}
\label{lem:trans_ident}For all $x,y,z$ in a Kikkawa left loop, 
\begin{equation}
L_{x}y\cdot (L_{x}z)^{-1}=L(x,y)(y\cdot z^{-1}).  \label{eq:trans_ident}
\end{equation}
\end{lemma}

\begin{proof}
In any LIP\ left loop, we have the identity 
\begin{equation}
L(x,y)^{-1}=L(x^{-1},x\cdot y)  \label{eq:LIP_inv}
\end{equation}
for all $x,y$ (see, e.g., \cite{kiechle}, 3.2). In any Kikkawa left loop, we
have the identities 
\begin{eqnarray}
L(x,y)^{-1} &=&L(y,x)  \label{eq:kik_inv} \\
L(x,y)(y\cdot x) &=&x\cdot y  \label{eq:bruck}
\end{eqnarray}
for all $x,y$ (see, e.g., \cite{kiechle}, 3.5(4) and 3.5(2)). We compute 
\begin{eqnarray*}
(x\cdot y)\cdot (x\cdot z)^{-1} &=&(x\cdot y)\cdot (x^{-1}\cdot z^{-1}) \\
&=&((x\cdot y)\cdot x^{-1})\cdot L(x\cdot y,x^{-1})z^{-1}.
\end{eqnarray*}
Now 
\begin{eqnarray*}
(x\cdot y)\cdot x^{-1} &=&L(x\cdot y,x^{-1})(x^{-1}\cdot (x\cdot y)) \\
&=&L(x\cdot y,x^{-1})y,
\end{eqnarray*}
using (\ref{eq:LIP_inv}) and LIP. Thus 
\begin{equation*}
(x\cdot y)\cdot (x\cdot z)^{-1}=L(x\cdot y,x^{-1})(y\cdot z^{-1}),
\end{equation*}
using the $A_{l}$ property. Finally by (\ref{eq:kik_inv}), (\ref{eq:LIP_inv}%
), and (\ref{eq:kik_inv}) again, 
\begin{equation*}
L(x\cdot y,x^{-1})=L(x^{-1},x\cdot y)^{-1}=L(x,y)^{-1}=L(y,x).
\end{equation*}
This establishes the result.
\end{proof}

\begin{theorem}
\label{thm:distance_props}The following hold.

\begin{enumerate}
\item  For all $\mathbf{x},\mathbf{y}\in \mathbb{S}$, $0\leq d_{s}(\mathbf{x}%
,\mathbf{y})\leq \pi $. In addition, $d_{s}(\mathbf{x},\mathbf{y})=0$ if and
only if $\mathbf{x}=\mathbf{y}$, and $d_{s}(\mathbf{x},\mathbf{y})=\pi $ if
and only if $\mathbf{x}=-\mathbf{y}$.

\item  (Triangle inequality) For all $\mathbf{x},\mathbf{y},\mathbf{z}\in 
\mathbb{S}$, $d_{s}(\mathbf{x},\mathbf{y})\leq d_{s}(\mathbf{x},\mathbf{z}%
)+d_{s}(\mathbf{z},\mathbf{y})$.

\item  (Invariance under $O(\mathbf{V})$) For all $\mathbf{x},\mathbf{y}\in 
\mathbb{S}$, $d_{s}(A\mathbf{x},A\mathbf{y})=d_{s}(\mathbf{x},\mathbf{y})$.

\item  (Invariance under $L(\mathbb{S})$) For all $\mathbf{x},\mathbf{y},%
\mathbf{z}\in \mathbb{S}$, $d_{s}(L_{\mathbf{x}}\mathbf{y},L_{\mathbf{x}}%
\mathbf{z})=d_{s}(\mathbf{y},\mathbf{z})$.
\end{enumerate}
\end{theorem}

\begin{proof}
(1) follows from Theorem \ref{thm:norm_props}(1). (3) follows from Theorems 
\ref{thm:automorphisms} and \ref{thm:norm_props}(4). For (2), we have 
\begin{equation*}
\mathbf{x}\odot \mathbf{y}^{-1}=\mathbf{x}\odot (\mathbf{z}^{-1}\odot (%
\mathbf{z}\odot \mathbf{y}^{-1}))=(\mathbf{x}\odot \mathbf{z}^{-1})\odot L(%
\mathbf{x},\mathbf{z}^{-1})(\mathbf{z}\odot \mathbf{y}^{-1}).
\end{equation*}
Taking norms of both sides, and then using Theorem \ref{thm:norm_props}%
(3,4), we obtain the desired result. For (4), we use (\ref{eq:trans_ident})
in $(\mathbb{S},\odot )$, take norms, and apply Theorem \ref{thm:norm_props}%
(4): 
\begin{eqnarray*}
d_{s}(\mathbf{x}\odot \mathbf{y},\mathbf{x}\odot \mathbf{z}) &=&\left\| (%
\mathbf{x}\odot \mathbf{y})\odot (\mathbf{x\odot z})^{-1}\right\| _{s} \\
&=&\left\| L(\mathbf{x},\mathbf{y})(\mathbf{y\odot z}^{-1})\right\| _{s} \\
&=&\left\| \mathbf{y\odot z}^{-1}\right\| _{s} \\
&=&d_{s}(\mathbf{y},\mathbf{z}).
\end{eqnarray*}
This completes the proof.
\end{proof}

Theorem \ref{thm:norm_props}(3,4) and the factorization $O(\mathbf{H})=L(%
\mathbb{S})O(\mathbf{V})$ show that the distance function $d_{s}:\mathbb{S}%
\times \mathbb{S}\rightarrow \lbrack 0,\infty )$ is invariant under the
action of the entire orthogonal group $O(\mathbf{H})$ on $\mathbb{S}$.

To conclude, we note that the $\mathbb{R}$-odule structure $(\mathbb{S}%
,\odot ,\omega )$ can be used to parametrize some of the interesting curves
of spherical geometry; cf. the remarks in \cite{ungar-axioms} about curves
in hyperbolic geometry. For instance, for $\mathbf{x},\mathbf{y}\in \mathbb{S%
}$ with $\mathbf{x}\neq -\mathbf{y}$, the unique spherical line (great
circle) through $\mathbf{x}$ and $\mathbf{y}$ turns out to be given by 
\begin{equation}
\gamma (t)=\mathbf{x}\odot (\mathbf{x}^{-1}\odot \mathbf{y})^{t},
\label{eq:sphere_line}
\end{equation}
$t\in \mathbb{R}$. (Note that $\gamma (0)=\mathbf{x}$ and $\gamma (1)=%
\mathbf{y}$ by LIP.) Another interesting curve through $\mathbf{x},\mathbf{y}%
\in \mathbb{S}$ where $\mathbf{x}\neq -\mathbf{y}^{-1}$ and $(\mathbf{y}%
\odot \mathbf{x})\odot \mathbf{y}^{-1}\neq -\mathbf{y}$ is given by 
\begin{equation}
\eta (t)=(\mathbf{y}\odot ((\mathbf{y}\odot \mathbf{x})^{-1}\odot \mathbf{y}%
))^{t}\odot \mathbf{x},  \label{eq:equi_curve}
\end{equation}
$t\in \mathbb{R}$. (Note that $\eta (0)=\mathbf{x}$ and $\eta (1)=\mathbf{y}$
by Theorem \ref{thm:bol}. The conditions on $\mathbf{x},\mathbf{y}$
guarantee that Theorem \ref{thm:bol} can be applied.) The curve $\eta $ is
also a circle on $\mathbb{S}$ which is everywhere equidistant from the
spherical line $\nu (t)=(\mathbf{y}\odot ((\mathbf{y}\odot \mathbf{x}%
)^{-1}\odot \mathbf{y}))^{t}$. Indeed, $d_{s}(\eta (t),\nu (t))=d_{s}(%
\mathbf{x},\mathbf{e}_{0})=\left\| \mathbf{x}\right\| _{s}$, by Theorem \ref
{thm:distance_props}(4).

\section{Complex Spheres}

To conclude this paper, we mention that many of the preceding ideas can be
extended to \emph{complex} spheres. We illustrate this with the complex $1$%
-sphere 
\begin{equation}
S_{\mathbb{C}}^{1}=\left\{ \mathbf{x}=(x_{0},x_{1})\in \mathbb{C}^{2}:\left|
x_{0}\right| ^{2}+\left| x_{1}\right| ^{2}=1\right\} .  \label{eq:SC1}
\end{equation}
For $\mathbf{x},\mathbf{y}\in S_{\mathbb{C}}^{1}$, define 
\begin{equation}
\mathbf{x}\odot \mathbf{y}=\left\{ 
\begin{array}{cc}
\left( \dfrac{x_{0}}{\bar{x}_{0}}(\bar{x}_{0}y_{0}-\bar{x}%
_{1}y_{1}),x_{0}y_{1}+x_{1}y_{0}\right) & \text{if }x_{0}\neq 0 \\ 
\left( -\bar{x}_{1}y_{1},x_{1}y_{0}\right) & \text{if }x_{0}=0
\end{array}
\right. .  \label{eq:op2}
\end{equation}
It turns out that $(S_{\mathbb{C}}^{1},\odot )$ is an $A_{l}$, LIP left
loop. However, $(S_{\mathbb{C}}^{1},\odot )$ does \emph{not} satisfy AIP. A
fuller exploration of $(S_{\mathbb{C}}^{1},\odot )$ will appear elsewhere.

\end{document}